\documentclass[a4paper]{amsart}

\usepackage{latexsym, amssymb, amsfonts, eucal}
\usepackage{amsthm, amsmath}

\usepackage{mathrsfs}
\usepackage{hyperref}
\usepackage{enumerate}
\usepackage{calc}

\input xy
\xyoption{all}
\CompileMatrices

\title{The Kre\u{i}n-Mil'man Theorem for Metric Spaces with a Convex
Bicombing}
\author{Theo B\"uhler}
\address{Department of Mathematics, ETH Z\"{u}rich, CH-8092 Z{\"u}rich, %
  Switzerland}
\email{theo@math.ethz.ch}
\subjclass[2000]{51F99; 52A30}
\date{\today}

\DeclareMathOperator{\bbR}{\mathbb{R}}
\DeclareMathOperator{\CAT}{CAT}

\newtheorem{Lem}{Lemma}
\newtheorem*{Theorem}{Theorem}

\theoremstyle{definition}

\newtheorem*{Exms}{Examples}
\newtheorem*{Definition}{Definition}
\newtheorem*{Acknowledgment}{Acknowledgment}

\begin{document}

\begin{abstract}
  We use bicombings on arcwise connected
  metric spaces to give definitions of convex sets and extremal points.
  These notions coincide with the customary ones
  in the classes of normed vector spaces and geodesic metric spaces which are
  convex in the usual sense. A rather straightforward modification
  of the standard proof of the Kre\u\i{}n-Mil'man Theorem yields the
  result that in a large class of metric spaces every compact convex set is
  the closed convex hull of its extremal points. The result appears to be new
  even for $\CAT{(0)}$-spaces.
\end{abstract}

\maketitle

\begin{Definition}
  Let $(X,d)$ be a metric space.
  A \emph{convex bicombing} on $X$ is a map
  $X \times X \to C([0,1],X)$, $(x,y) \mapsto [x,y](\,\cdot\,)$
  satisfying
  \begin{enumerate}[(i)]
    \item
      For all $x,y \in X$, one has $[x,y](0) = x$ and $[x,y](1) = y$.
      Moreover, $[x,x] \equiv x$.
 
    \item
      For all $x,y, x', y' \in X$ the function
      $t \mapsto d \left([x,y](t), [x',y'](t)\right)$ is a convex function
      on $[0,1]$.
  \end{enumerate}
 % Observe that a metric space admitting a convex bicombing must be
 % path-connected.
\end{Definition}

\begin{Exms}
  We are interested in the following two special cases:
  \begin{enumerate}[(i)]
    \item
      Let $(X,d)$ be a \emph{convex} metric space in the sense that any
      two points of $X$ are connected by a geodesic and that the inequality
      \[
      d(c(t),c'(t)) \leq (1-t) d(c(0),c'(0)) + t d(c(1),c'(1))
      \]
      holds for all linearly parameterized geodesics $c,c':[0,1] \to X$
      and all $t \in [0,1]$.
      It is easy to see that such a space is uniquely geodesic
      and that the map associating to $x,y \in X$ the unique linearly
      parameterized geodesic $[x,y]$ connecting $x$ to $y$ is a
      convex bicombing. Observe that the class of convex metric spaces
      contains in particular all $\CAT{(0)}$-spaces.

    \item
      If $(X,\|\,\cdot\,\|)$ is a normed vector space then by setting
      $[x,y](t) = (1-t)x + ty$ for all $x,y \in X$ one
      obtains a convex bicombing on $X$.
  \end{enumerate}
\end{Exms}

Throughout this note, $X$ stands as a shorthand for a metric space
$(X,d,[\,\cdot\,,\,\cdot\,])$ with a convex bicombing.  The choice of a
bicombing in a metric space yields notions of convexity and extremity as
follows: A subset $C \subset X$ is called \emph{convex} if $x,y \in C$
implies $[x,y] \subset C$. A function $\phi: A \to \bbR$ on a subset $A$ of
$X$ is called \emph{convex} if $t \mapsto \phi([x,y](t))$ is a convex 
function on the interval $[0,1]$ for all $x,y \in A$ with $[x,y] \subset A$.
The \emph{(closed) convex hull} of a subset $A$ of $X$ is the smallest
(closed) convex set containing $A$.

\begin{Definition}
  Let $C \subset X$ be a convex set. A \emph{closed} non-empty subset
  $E \subset C$ is called
  \emph{extremal} if for all $x,y \in C$ for which there exists some
  $t \in (0,1)$ with $[x,y](t) \in E$ one has $[x,y] \subset E$.
  A point $p \in C$ is called extremal if $\{p\}$ is an extremal set.
\end{Definition}

We can now state the Kre\u{\i}n-Mil'man Theorem:

\begin{Theorem}
  If $C$ is a compact and convex subset of a metric space $X$ with
  convex bicombing then $C$ is the closed convex hull of the
  set of its extremal points.
\end{Theorem}

Consider the family of extremal subsets of $C$ and order it by 
inclusion. This family is non-empty since it contains $C$
and by compactness of $C$ the
intersection of a decreasing family of extremal sets is non-empty,
hence it is an extremal set. Zorn's lemma applies and yields:

\begin{Lem}
  Every extremal subset of $C$ contains a minimal extremal subset.\qed
\end{Lem}

The next step is to prove:

\begin{Lem}
  Every minimal extremal subset is a one-point set.
\end{Lem}

Before giving the proof, we record the following:

\begin{Lem}
  Let $E$ be an extremal subset of $C$ and let $\phi:E \to \bbR$ be
  continuous and convex. Then $E_\phi = \{e \in E\,:\,\phi(e) = \max{\phi}\}$
  is an extremal subset of $C$.
\end{Lem}

\begin{proof}[Proof of Lemma~3]
  Since $E$ is compact and $\phi$ is continuous,
  the set $E_\phi$ is closed and non-empty. Let $x,y \in C$ and
  $t \in (0,1)$ be such that $[x,y](t) \in E_\phi$. Since $E$ is extremal,
  we have $[x,y] \subset E$. Now notice that $t \mapsto \phi([x,y](t))$
  is convex on $[0,1]$ (by convexity of the bicombing)
  and assumes its maximum at some point $t \in (0,1)$,
  hence it is constant and thus $[x,y] \subset E_\phi$.
\end{proof}

\begin{proof}[Proof of Lemma~2]
  Let $E$ be a minimal extremal subset of $C$ and let $e \in E$.  Consider
  the function $\phi:x \mapsto d(x,e)$ which is convex by the convexity of the
  bicombing and the assumption $[x,x] \equiv x$.
  By Lemma~3 the set $E_\phi$ is an
  extremal subset of $C$. If there exists a point $e' \in E$ distinct from
  $e$, then $\phi$ is not constant, hence $E_\phi$ is a proper extremal
  subset of $E$ contradicting the minimality of $E$.
\end{proof}

We need one more fact before we can finish the proof of the Kre\u{\i}n-Mil'man
Theorem.

\begin{Lem}
  For a compact and convex set $K \subset X$ put
  $d_K (x) = \min_{k \in K} d(k,x)$. The function $d_K : X \to \bbR$ is
  continuous and convex.
\end{Lem}

\begin{proof}
  It follows from the triangle inequality that $d_K$ is $1$-Lipschitz, hence
  continuous. Let $x,y \in X$ be arbitrary points and pick
  points $\bar{x},\bar{y} \in K$ such that $d(x,\bar{x}) = d_K (x)$
  and $d(y,\bar{y}) = d_K (y)$. By convexity of the bicombing the function
  $t \mapsto d([x,y](t), [\bar{x},\bar{y}](t))$ is convex and by convexity
  of $K$, we have $[\bar{x},\bar{y}] \subset K$, so
  \[
  d_K ([x,y](t)) \leq d([x,y](t),[\bar{x},\bar{y}](t)) \leq
  (1-t) d(x,\bar{x}) + t d(y,\bar{y}) = (1-t) d_K (x) + t d_K (y)
  \]
  and it follows that $d_K$ is convex.
\end{proof}

\begin{proof}[Proof of the Theorem]
  Let $K$ be the closed convex hull of the set of extremal points of $C$.
  By Lemma~1 and Lemma~2, $K$ is non-empty and by convexity of $C$
  we have $K \subset C$. By Lemma~4 the function $\phi = d_K$ is continuous and
  convex on $C$, so by Lemma~3 $E_\phi$ is an extremal set, and it
  contains an extremal point of $C$ by Lemma~1 and Lemma~2,
  so $E_\phi \cap K$ is non-empty.
  If there existed a point $p \in C \smallsetminus K$ then
  $\phi$ would be non-constant, so $E_\phi$ would have to be disjoint from $K$,
  a contradiction.
\end{proof}

\begin{Acknowledgment}
  I would like to thank Stefan Wenger for reminding me of the many interesting
  problems concerning convex hulls in $\CAT{(0)}$-spaces. Arvin Moezzi and
  Stefan Wenger pointed out a slip in the definition of a convex bicombing
  contained in an earlier version of this paper. 
\end{Acknowledgment}

\end{document}